\documentclass[12pt]{article}%
\usepackage{amssymb}
\usepackage{amsfonts}
\usepackage{amsmath}
\usepackage{graphicx}%
\begin{document}

\begin{center}
{\Large The Kraft sum as a monotone function on the refinement-ordered set }

{\Large of uniquely decipherable codes }

\bigskip

Stephan Foldes

Tampere University of Technology

\end{center}

Abstract. \textit{The set of all} \textit{uniquely decipherable (UD) codes is
partially ordered by refinement, meaning that all strings in the cruder code
can be represented as \ \ concatenations of strings taken from the finer code.
The Kraft sum is a monotone (increasing) function on this poset. } \textit{In
the refinement order, chains of UD codes having the same Kraft sum are
necessarily of the simple descending type.}

\bigskip

\bigskip

\textbf{Introduction}

\bigskip

Let $A$ be any non-empty finite set, called \textit{alphabet}. Let
$con:A^{\ast\ast}\longrightarrow A^{\ast}$ be the concatenation map which to
every string of strings associates their concatenation. In this note we shall
call any finite subset $C\subseteq A^{\ast}$ not containing the null string a
\textit{code. }A code is said to be \textit{uniquely decipherable (UD) }if
$con$ is injective on the subset $C^{\ast}$ of $A^{\ast\ast}$. (Note that
several authors, including Berstel and Perrin [1], reserve the term "code" to
mean UD code.)

For any codes $C$ and $D$ write $C$ $\leq$ $D$ and say that $D$ is
\textit{finer} than $C$, or that it is a \textit{refinement} of $C,$ if
$C\subseteq con[D].$ This is a partial order relation (antisymemtry is ensured
by unique decipherability). We say that $D$ is an \textit{irredundant
refinement} of $C$ if $C\leq D$ and no proper subset of $D$ is finer than $C$.
Every code has infinitely many refinements. However, due to its finiteness,
each code can have only finitely many irredundant refinements.

Denoting by $r$ the number of elements of the \textit{alphabet} $A$, the
\textit{Kraft sum} $K(C)$ of any code $C\subseteq A^{\ast}$ is defined as
$\sum_{\mathbf{x}\in C}$ $r^{-len(\mathbf{x})},$ where $len$ is the length
function. In [5] McMillan showed that if $C$ is a uniquely decipherable code,
then its Kraft sum is at most $1.$ Simplified combinatorial proofs were given
by Karush [3], and by Berstel and Perrin ([1] Chapter 1, Theorem 4.2). The
proof was also reformulated in [2] as an argument involving evaluations of
polynomials with non-commuting indeterminates corresponding to the various
(infinitely many) strings in $A^{\ast}.$ In [2] we also concluded that for any
UD codes $C$ and $D$ such that $C\leq D,$ the inequality $K(C)\leq K(D)$
holds. Here the purely combinatorial proof due to Berstel and Perrin [1] is
shown to yield the same conclusion that was reached in [2], and from this some
further conclusions are drawn about the set of UD codes having the same Kraft sum.

The Kraft inequality as originally established by Kraft [4] stated that
$K(C)\leq1$ for instantaneous (prefix-free) codes, which are a special class
of UD codes. This can be verified in several ways - for a recent approach,
which also generalizes the inequality to data structures other than strings,
see Valmari [6].

\bigskip

\bigskip

\textbf{Statements and proofs}

\textbf{\bigskip}

For any code $C$ and positive integer $k$ denote by $C^{k}$ the code
consisting of all possible concatenations of $k$ (not necessarily distinct)
members of $C.$ Note that if $C$ is a UD code, then $C^{k}$ is also UD and
$Card(C^{k})=[Card(C)]^{k}.$ The following appears in Berstel and Perrin [1],
Chapter 1, Proposition 4.1.

\bigskip

\textbf{Proposititon 1 } (from [1])\ \ \textit{For any code} $C$ \textit{over
a given alphabet} \textit{and positive integer} $k,$\ \textit{we have}
$K(C^{k})\leq K(C)^{k}$ \textit{and\ the following conditions are equivalent:}

(i) $C$\ \textit{is uniquely decipherable,}

(ii) $K(C^{k})=K(C)^{k}$ \textit{for all positive integers} $k.$

\textit{Proof.} \ (As given by Berstel and Perrin [1].)\ For each positive
$k$, denote by $C^{(k)}$ the set of strings $(\mathbf{v}_{1},...,\mathbf{v}%
_{k})$ of $k$ (not necessarily distinct) words from $C.$Clearly,
$C^{(k)}\subseteq A^{\ast\ast}.$

We claim that $C$ is UD if and only if concatenation restricted to $C^{(k)}$
is an injective map for every $k.$ These injectivity conditions are clearly
necessary for $C$ to be UD. On the other hand, if $C$ is not UD, then for some
positive $m,n$ two different strings of words, $x=(\mathbf{x}_{1}%
\mathbf{,...,x}_{m})\in C^{(m)}$ and $y=(\mathbf{y}_{1}\mathbf{,...,y}_{n})\in
C^{(n)}$ yield the same concatenation, $con$ $x=con$ $y.$ Let $k=m+n$. The
strings of words $(\mathbf{x}_{1}\mathbf{,...,x}_{m},\mathbf{y}_{1}%
\mathbf{,...,y}_{n})$ and $(\mathbf{y}_{1}\mathbf{,...,y}_{n},\mathbf{x}%
_{1}\mathbf{,...,x}_{m})$ are both in $C^{(k})$, they are distinct, and they
yield the same concatenation, proving the claim.

Observe that concatenation restricted to $C^{(k)}$ is always a surjective map
onto $C^{k}.$

For $(\mathbf{v}_{1},...,\mathbf{v}_{k})\in C^{(k)}$ the word $con(\mathbf{v}%
_{1},...,\mathbf{v}_{k})$ in $C^{k}$ contributes to the Kraft sum $K(C^{k})$ a
term equal to the product of the terms $r^{-len(\mathbf{v}_{i})}$, $1\leq
i\leq k$ of the Kraft sum of $C.$ Adding up these products over all members of
$C^{(k)}$ equals $K(C)^{k}$\ and it yields exactly $K(C^{k})$ if the map $con$
is injective on $C^{(k)}$, otherwise it yields a strict upper bound of
$K(C^{k})$.
$\ \ \ \ \ \ \ \ \ \ \ \ \ \ \ \ \ \ \ \ \ \ \ \ \ \ \ \ \ \ \ \ \square$

\bigskip

We can apply to two codes comparable by refinement the reasoning presented in
Berstel and Perrin's proof of Theorem 4.2 in Chapter 1 of [1].\ Let\ $C$ be a
UD code and let $D$ be any code finer than $C$. \ There is a positive integer
$m$ such that for all integers $n>m,$ $C$ is disjoint from $D^{n}.$ For any
fixed positive integer $k$, we have (by an obvious induction with respect to
$k$) that%
\[
C^{k}\subseteq D^{k}\cup D^{k+1}\cup...\cup D^{mk}%
\]
and $K(C)^{k}=K(C^{k})$ is less than or equal to
\begin{align*}
&  K(D^{k}\cup D^{k+1}\cup...\cup D^{mk})\\
&  \leq K(D^{k})+K(D^{k+1})+...+K(D^{mk}\\
&  \leq K(D)^{k}+K(D)^{k+1}+...+K(D)^{mk}\\
&  =K(D)^{k}[1+K(D)+...+K(D)^{(m-1)k}]
\end{align*}
First, when we choose the finest code $D$ consisting of all the words of
length $1$, we get $K(C)^{k}\leq(m-1)k+1$. It follows that $K(C)\leq\left[
(m-1)k+1\right]  ^{1/k}$, for every $k.$ Necessarily, $K(C)\leq1.$ Second, if
we assume that code $D$ is also UD, and thus it also has Kraft sum at most
$1$, we get%
\[
\left[  \frac{K(C)}{K(D)}\right]  ^{k}\leq1+.K(D)+...+K(D)^{(m-1)k}%
\leq(m-1)k+1
\]
which is true for all $k$, implying that the ratio on the left hand side is at
most $1$, $K(C)\leq K(D).$ This yields the following extension of McMillan's Theorem:

\bigskip

\textbf{Proposition 2} \textit{\ The Kraft sum is a monotone (increasing)
function on the refinement-ordered set of uniquely decipherable codes. For
each UD code }$C$ \textit{there are only finitely many finer UD codes with the
same Kraft sum. \ \ \ }$\ \ \square$

\bigskip\ 

The second statement follows from the fact, noted above, that a UD code $C$
has only finitely many irredundant refinements, and from the observation that
for a code $D$ finer than $C$, and which is not an irredundant refinement of
$C$, the Kraft sums cannot be equal, the inequality between them has to be
strict, $K(C)<K(D).$ By Proposition 2, every UD code $C$ has at least one UD
refinement $D$ with the same Kraft sum, and which can no longer be properly
refined without increasing the Kraft sum. As a further consequence, we have:

\bigskip

\textbf{Proposition 3 \ }\textit{In the refinement-ordered set of UD codes,
all infinite chains of UD codes with the same Kraft sum are of type }%
$\omega^{\ast}$ \textit{(i.e. of the same order type as the negative
integers). \ \ }$\ \ \square$

\bigskip

Infinite chains of UD codes all having the same Kraft sum exist indeed, in
fact there is such a chain below every member of the poset of UD codes: for
any UD code $C$, consider for example $C>C^{2}>C^{4}>...>C^{2^{n}}>...$

\bigskip

\textbf{Acknowledgements. }

This work has been co-funded by Marie Curie Actions and supported by the
National Development Agency (NDA) of Hungary and the Hungarian Scientific
Research Fund (OTKA, contract number 84593), within a project hosted by the
University of Miskolc, Department of Analysis.

The author wishes to thank S. Fegyverneki, S. Radeleczki, J. Szigeti and A.
Valmari for useful comments and discussions.

\bigskip

\includegraphics[height=15mm, width=20mm]{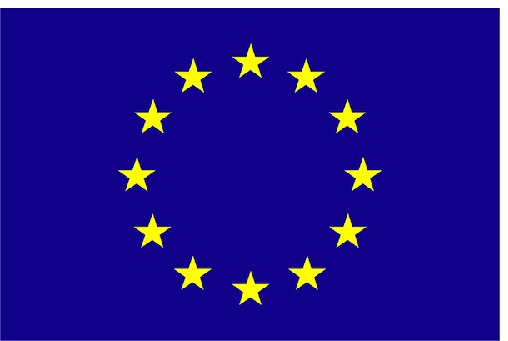}
\includegraphics[height=15mm, width=20mm]{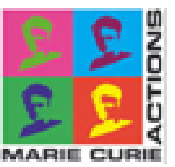}
\includegraphics[height=15mm, width=20mm]{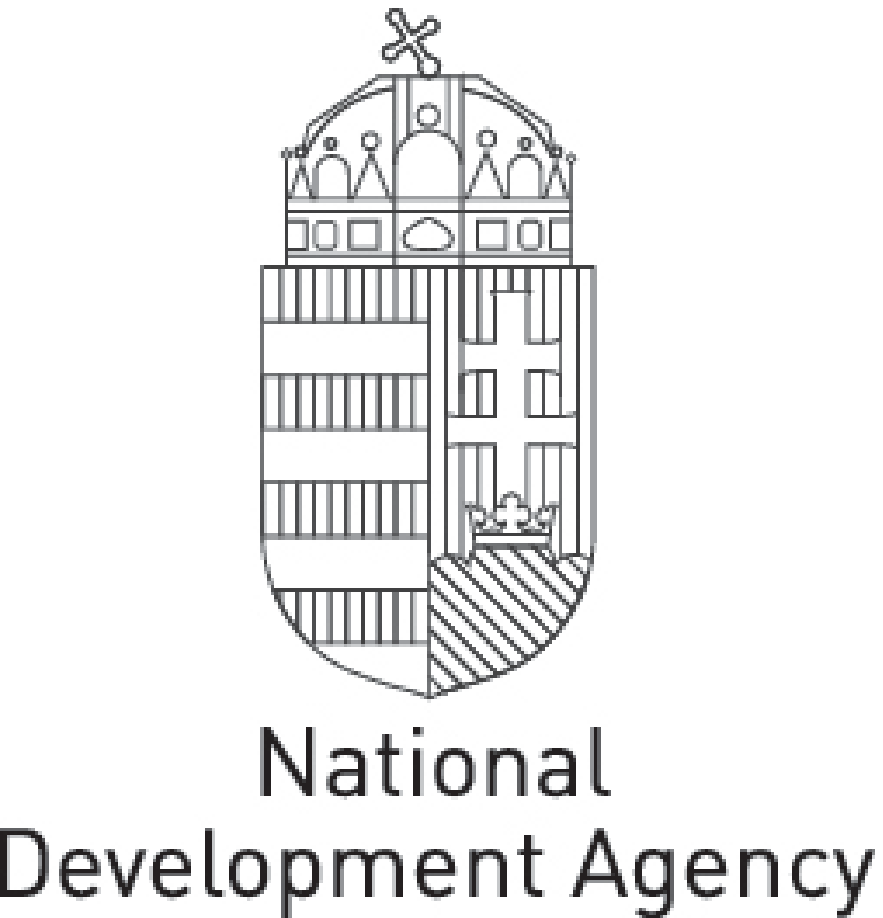} \includegraphics[height=15mm, width=20mm]{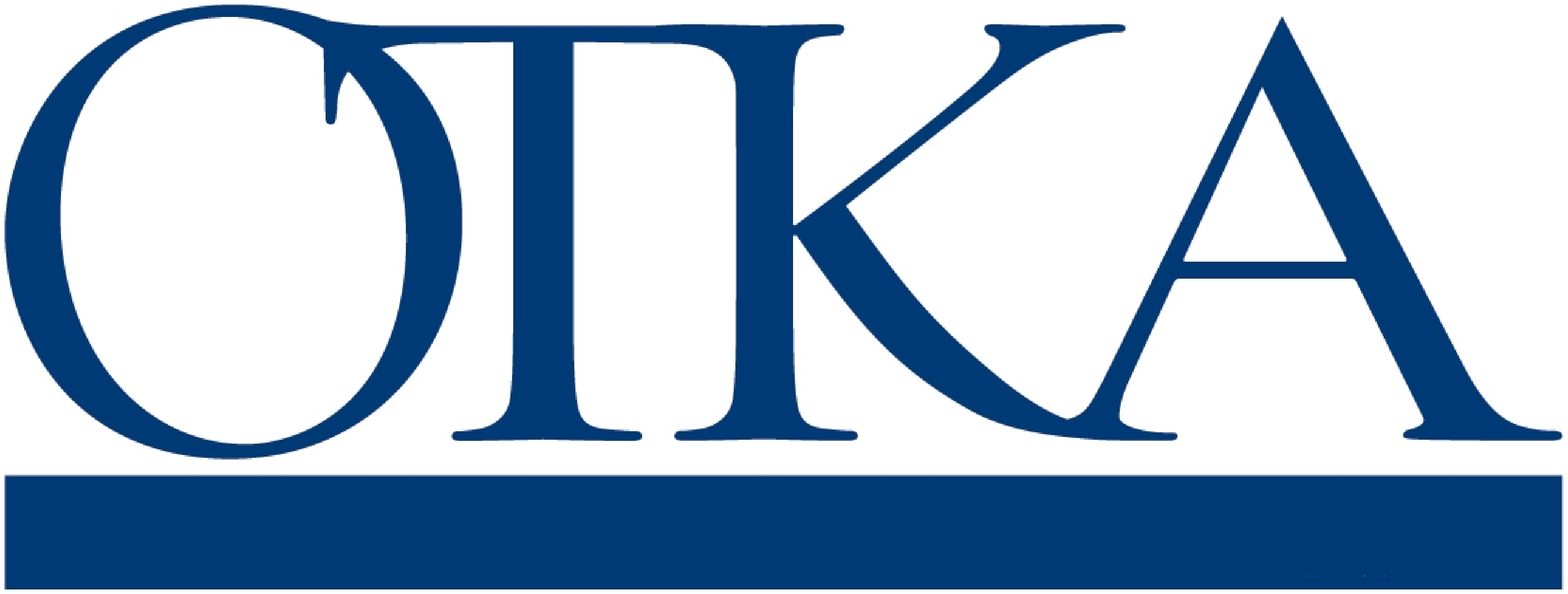}

\bigskip

\textbf{References}

\bigskip

[1] \ J. Berstel and D. Perrin, \textit{Theory of Codes}, Academic Press, 1985.

[2] \ S. Foldes, On McMillan's theorem about uniquely decipherable codes,
ArXiv:0806.3277v2 (2008)

[3] \ J. Karush, A simple proof of an inequality of McMillan, \textit{IRE
Trans. Information\ Theory}\textbf{ }IT-7 (1961), 118-118

[4] \ L.G. Kraft, \textit{A Device for Quantizing, Grouping, and Coding
Amplitude Modulated Pulses}, Q.S. Thesis, MIT, 1949

[5] \ B. McMillan, Two inequalities implied by unique decipherability,
\textit{IRE Trans. Information Theory} IT-2 (1956), 115-116

[6] \ A. Valmari, Does the Shannon bound really apply to data structures?
\textit{Proc. Estonian Acad. Sc.} 62, No. 1 (2013), 47-58

\bigskip

\end{document}